\documentclass[11pt]{article}

\usepackage{a4wide}
\usepackage{latexsym}
\usepackage{amsfonts}
\usepackage{amsthm}
\usepackage{amsmath}
\usepackage{amsopn,amscd, amssymb}

\def\HH{\mathrm H}
\def\piN{\pi_N}
\def\eG{\mathrm e^M_G}
\def\eMN{\mathrm e}
\def\EN{E_N}

\newtheorem{Theorem}{Theorem}
\newtheorem{Remark}[Theorem]{Remark}
\newtheorem{thm}{Theorem}[section]
\newtheorem{prop}[thm]{Proposition}
\newcommand{\Aut}{\mathop{\mbox{\rm Aut}}\nolimits}
\newcommand{\Out}{\mathop{\mbox{\rm Out}}\nolimits}
\newcommand{\Der}{\mathop{\mbox{\rm Der}}\nolimits}

\numberwithin{equation}{section}
\title
{Exact sequences in the cohomology of a group extension}
\author{
J.~Huebschmann
\\[0.3cm]
 USTL, UFR de Math\'ematiques\\
CNRS-UMR 8524
\\
Labex CEMPI (ANR-11-LABX-0007-01)
\\
59655 Villeneuve d'Ascq Cedex, France\\
Johannes.Huebschmann@math.univ-lille1.fr
 }

\long
\def\MSC#1\EndMSC{\def\arg{#1}\ifx\arg\empty\relax\else
      {\par\narrower\noindent
      2010 Mathematics Subject Classification: #1\par}\fi}
\long
\def\KEY#1\EndKEY{\def\arg{#1}\ifx\arg\empty\relax\else
    {\par\narrower\noindent
      Keywords and Phrases: #1\par}\fi}

\begin{document}

\maketitle
\medskip
\centerline
{Dedicated to the memory of my friend and  colleague 
Chih-Han Sah}
\medskip

\begin{abstract}
In \cite{MR2959787}
the authors constructed a seven term exact sequence in the 
cohomology of a group extension
$1 \to N \to G \to Q \to 1$ with coefficients in a $G$-module
$M$. However, they were unable 
to establish the precise
link between the maps in that sequence and the 
corresponding maps arising from the spectral
sequence associated to the group extension and the $G$-module $M$. 
In this paper, we show that there is a close connection
between  \cite{MR2959787} and the two earlier papers
\cite{MR641328} and
\cite{MR597986}. In particular, we show that 
the results in \cite{MR641328} and
\cite{MR597986} entail that
the maps of 
\cite{MR2959787}
other than the obvious inflation and restriction
maps do correspond to the corresponding ones arising from
the spectral sequence.

\end{abstract}

\MSC 18G40 20J05 20J06 55R20
\EndMSC

\KEY Group extension, 
group cohomology, 
low degree exact sequences, spectral sequence of a group extension
\EndKEY

\tableofcontents

\section{Introduction}
Consider an extension of (discrete) groups
\begin{equation}
\begin{CD}
\mathrm{e}_G\colon
1
@>>>
N
@>{i^N}>>
G
@>{\pi}>>
Q
@>>>
1
\end{CD}
\label{ge1}
\end{equation}
and a $G$-module $M$.
Let $(\mathrm E^{p,q}_r(M),d_r)$
denote the spectral sequence associated to the group extension
\eqref{ge1} and the $G$-module $M$.
As usual, for $K=N$ as well as for $K=Q$,
we will use the notation
 ${-}^K$ for the invariants relative to the group $K$,
and 
the notation \lq\lq $\mathrm{res}$\rq\rq\  and \lq\lq $\mathrm{inf}$\rq\rq\ 
will
refer to the maps in cohomology induced by the injection of $N$ into $G$
and by the projection from $G$ to $Q$, respectively.
Let
$\HH^2(G,M)_1$ denote the kernel of 
$\mathrm{res}\colon \HH^2(G,M)
\to
\HH^2(N,M)^Q$.
It is a straightforward exercise to show that the canonical
homomorphism 
\[
\mathrm{can}\colon
\HH^2(G,M)_1
\longrightarrow
\HH^1(Q,\HH^1(N,M))\cong \mathrm E_2^{1,1}(M)
\]
induced by the canonical surjection
$\HH^1(Q,\HH^1(N,M)) \to  \mathrm E_{\infty}^{1,1}(M)$,
cf. \eqref{can} below,
extends the classical five term exact sequence
in the cohomology of the group extension
\eqref{ge1}
with coefficients in $M$
by two more terms to 
an exact sequence of the kind
\begin{equation}
\begin{CD}
0
@>>>
\HH^1(Q,M^N)
@>{\mathrm{inf}}>>
\HH^1(G,M)
@>{\mathrm{res}}>>
\HH^1(N,M)^Q
\\
@>{\mathrm{d_2}}>>
\HH^2(Q,M^N)
@>{\mathrm{inf}}>>
\HH^2(G,M)_1
@>{\mathrm{can}}>>
\HH^1(Q,\HH^1(N,M))
\\
@>{\mathrm{d_2}}>>
\HH^3(Q,M^N).
\end{CD}
\label{exact1}
\end{equation}
In \cite[(3) p.~257]{MR0393273},
C. H. Sah spelled out this sequence,
with the notation \lq\lq tg\rq\rq\ 
rather than \lq\lq can\rq\rq.

Let
\begin{equation}
\lambda\colon \HH^1(Q,\HH^1(N,M)) \longrightarrow
\HH^3(Q,M^N)
\end{equation}
be the map constructed in \cite{MR641328}
that yields an explicit description of
the differential
\begin{equation*}
\mathrm{d_2}\colon 
\HH^1(Q,\HH^1(N,M)) \longrightarrow 
\HH^3(Q,M^N).
\end{equation*}
In \cite{MR2959787},
the authors constructed explicit maps
\begin{align}
\mathrm{tr}&\colon 
\HH^1(N,M)^Q \longrightarrow 
\HH^2(Q,M^N)
\label{tr}
\\
\rho&\colon 
\HH^2(G,M)_1 \longrightarrow 
\HH^1(Q,\HH^1(N,M))
\label{rho}
\end{align}
such that, with $\mathrm{tr}$, $\rho$, and $\lambda$
substituted for, respectively,
\begin{align}
\mathrm{d_2}&\colon 
\HH^1(N,M)^Q \longrightarrow 
\HH^2(Q,M^N),
\label{d21}
\\
\mathrm{can}&\colon 
\HH^2(G,M)_1 \longrightarrow 
\HH^1(Q,\HH^1(N,M)),
\\
\mathrm{d_2}&\colon 
\HH^1(Q,\HH^1(N,M)) \longrightarrow 
\HH^3(Q,M^N),
\label{d22}
\end{align}
an exact sequence of the kind
\eqref{exact1} results.
They noted \cite[p.~71, 4th paragraph]{MR2959787} 
that they did not know whether
 the maps
$\mathrm{tr}$ and $\rho$ 
in their sequence coincided with the 
maps 
$\mathrm{d_2}\colon 
\HH^1(N,M)^Q \to 
\HH^2(Q,M^N)$ and
$\mathrm{can}\colon 
\HH^2(G,M)_1 \to 
\HH^1(Q,\HH^1(N,M))$,
respectively, that are
obtained from the spectral sequence. 
In the addendum
\cite{MR2995036} to \cite{MR2959787},
they compared
the conceptual description of \eqref{tr}
given in 
their paper 
with a conceptual description of a map of the same kind
given in
\cite{MR597986} (Section 1.5)
and concluded that the descriptions obtained
were different.

In this note we shall show that 
the maps \eqref{tr} and \eqref{rho}
are special cases of maps given  in \cite{MR641328} and
\cite{MR597986}, see Subsection \ref{section5} 
and Section \ref{d23} below,
that
the requisite explicit descriptions
may be found in \cite{MR641328} and
\cite{MR597986}, and that
the results in
\cite{MR641328} and
\cite{MR597986} entail
the following.

\begin{Theorem}
\label{one}
The homomorphism 
{\rm \eqref{tr}}
constructed in {\rm\cite{MR2959787}}
coincides with  the
differential 
$\mathrm{d_2}\colon 
\HH^1(N,M)^Q \to 
\HH^2(Q,M^N)$ 
of the spectral sequence $(\mathrm E^{p,q}_r(M),d_r)$.
\end{Theorem}

\begin{Theorem}
\label{two}
The homomorphism  
{\rm\eqref{rho}}
constructed in {\rm \cite{MR2959787}}
coincides with the 
homomorphism
$\mathrm{can}\colon 
\HH^2(G,M)_1 \to 
\HH^1(Q,\HH^1(N,M))$
arising from the spectral sequence  $(\mathrm E^{p,q}_r(M),d_r)$.
\end{Theorem}

This settles
the questions 
left open in \cite{MR2959787} and
\cite{MR2995036}.
It also shows that the claim in \cite[p.~71]{MR2959787}
to the effect that, in the exact sequence constructed in that paper,
 there were no explicit descriptions available for the maps
other than the inflation map, restriction map and the map $\lambda$
is unfounded.
In particular, \cite[Proposition 3.1]{MR641328}
entails an explicit conceptual description
of the differential \eqref{d21},
and another conceptual description
of this differential
is given in \cite[Theorem 6, Section 5]{MR597986}.
Moreover, in Subsection \ref{compar} below we shall show that
the description of 
\eqref{tr}
given in \cite{MR2959787} is essentially the same as
that in our paper \cite{MR597986}.
For intelligibility we note that 
the main result in
\cite{MR597986} is an extension of the classical five term exact sequence
by three more terms, viz.
\begin{equation}
\begin{CD}
@. @. \ldots
@>>>
\HH^1(N,M)^Q
\\
@>{\Delta}>>
\HH^2(Q,M^N)
@>{\mathrm{inf}}>>
\HH^2(G,M)
@>{\mathrm{j}}>>
\mathrm{Xpext}(G,N;M)
\\
@>{\Delta}>>
\HH^3(Q,M^N)
@>{\mathrm{inf}}>>
\HH^3(G,M).
\end{CD}
\label{exact3}
\end{equation}
The description of \eqref{rho} given in \cite{MR2959787}
can be deduced from the construction of the map $j$, see
Section \ref{d23} below.
The conceptual descriptions of the differentials in
\cite{MR641328} and the
construction of the exact sequence
\eqref{exact3}
rely on the interpretation of group cohomology in terms of
crossed modules developed in \cite{crossed}.

We dedicate this paper to the memory of Chih-Han Sah, 
see, e.g., \cite{MR1490540}. 
Sah contributed substantially to the understanding of
the spectral sequence of a group extension, cf. \cite{MR0393273},
\cite {MR0463264}.

\section{$\mathrm{d_2}\colon 
\HH^1(N,M)^Q \longrightarrow 
\HH^2(Q,M^N)$}

We adjust the notation to 
the present discussion.
We denote the group ring of $G$ by $\mathbb Z G$,
the standard augmentation map by $\varepsilon \colon \mathbb ZG \to \mathbb Z$,
 and
the augmentation ideal $\ker(\varepsilon)$ 
by $\mathrm I G$, and we use the same kind of 
notation for the group ring of $Q$ and the augmentation ideal of $Q$.

\subsection{Description in \cite{MR641328}}

Let $T=\ker(\mathrm{Hom}_N(IG,M) \to \mathrm{Ext}_N^1(\mathbb Z,M)\cong \HH^1(N,M))$; with respect to the obvious $Q$-module structures on
$\mathrm{Hom}_N(IG,M)$ and $\mathrm{Ext}_N^1(\mathbb Z,M)$,
the map from $\mathrm{Hom}_N(IG,M)$ to $\mathrm{Ext}_N^1(\mathbb Z,M)$
is a morphism of $Q$-modules, and hence
$T$ inherits a $Q$-module structure.

Consider an extension
\begin{equation}
\begin{CD}
\mathrm{e}_M\colon
0
@>>>
M
@>{i_M}>>
E_M
@>{p_M}>>
\mathbb Z
@>>>
0
\end{CD}
\label{extM}
\end{equation}
of $N$-modules 
that
represents a member of $\mathrm{Ext}_{N}(\mathbb Z,M)^Q\cong
\HH^1(N,M)^Q$.
Lift the identity of $\mathbb Z$ to a commutative diagram
\begin{equation}
\begin{CD}
0
@>>>
\mathrm I G
@>>>
\mathbb Z G
@>{\varepsilon}>>
\mathbb Z
@>>>
0
\\
@.
@V{\mu}VV
@V{\nu}VV
@|
@.
\\
0
@>>>
M
@>{i_M}>>
E_M
@>{p_M}>>
\mathbb Z
@>>>
0
\end{CD}
\label{diag11}
\end{equation}
in the category of $N$-modules.
For $x \in G$, let $\alpha_x$ denote the associated operator on
$\mathrm IG$ as well as that on $M$ where the notation is slightly abused;
given $q \in Q$, pick $x \in G$ such that $\pi(x)=q$, and
define the derivation 
\begin{equation}
d \colon Q \longrightarrow \mathrm{Hom}_N(\mathrm IG,M),
\ 
d(q)=\alpha_x \mu \alpha_x^{-1}-\mu, \ q\in Q.
\label{T}
\end{equation}
Each value $d(q)$, as $q$ ranges over $Q$, is well defined, that is,
does not depend on the choice of $x\in G$ such that $\pi(x)=q$, and
the values of $d$ lie in $T$ since $[\mathrm{e}_M]\in 
\mathrm{Ext}_{N}(\mathbb Z,M)$ is $Q$-invariant.
Consider the semi-direct fiber product group
\begin{equation}
\mathrm{Hom}_N(\mathbb Z G,M)\rtimes _TQ
=\{(\varphi,q);\varphi|_{\mathrm I G}=d(q)\}\subseteq 
\mathrm{Hom}_N(\mathbb Z G,M)\rtimes Q .
\label{sdfb}
\end{equation}
The projection to $Q$ yields the group extension
\begin{equation}
\begin{CD}
\widehat {\mathrm{e}}_M\colon 
0
@>>>
M^N
@>>>
\mathrm{Hom}_N(\mathbb Z G,M)\rtimes _TQ
@>>>
Q
@>>>
1.
\end{CD}
\label{extQ}
\end{equation}
In view of 
\cite[Proposition 3.1]{MR641328},
the assignment to $\mathrm{e}_M$ of the class 
$[\widehat {\mathrm{e}}_M]\in \HH^2(Q,M^N)$
represented by $\widehat {\mathrm{e}}_M$ yields a conceptual description of
the differential \eqref{d21}.

\begin{Remark}
{\rm
A cocycle description  of
the differential \eqref{d21} is given in 
\cite[Remark 3.2]{MR641328}.
For the special case where $N$ acts trivially on $M$,
a similar cocycle description can be found in 
\cite[p.~21]{MR0463264}. In the latter paper, C. H. Sah explored 
some of the differentials in the spectral sequence
of a split group extension via certain characteristic classes.
These characteristic classes were extended in \cite{MR994078};
in that paper, the interested reader will also find more references
related with these characteristic classes.}
\end{Remark}

\subsection{Description in \cite{MR597986}} 
\label{hueb}

We now recall the construction in \cite[Section 5]{MR597986}.
Let $\overline{\mathrm{Aut}(\mathrm{e}_M)}$ denote the group of pairs
$(\alpha,x)$ such that $\alpha$ is an automorphism of $E_M$
whose restriction $\alpha_M$
to $M$ is an automorphism of $M$ and coincides with
the automorphism $x_M$ of $M$ induced by $x \in G$
via the $G$-action on $M$. Since the class 
$[\mathrm{e}_M]\in \mathrm{Ext}_{N}(\mathbb Z,M)$ 
is $Q$-invariant, 
every $x\in G$ lifts to a member $(\alpha,x)$ of 
$\overline{\mathrm{Aut}(\mathrm{e}_M)}$.
Moreover, given $m \in M$, let 
$\alpha_m\colon E_M \to E_M$ be the automorphism of $E_M$
given by
\begin{equation*}
\alpha_m(y)=y+ p_M(y)m,\ y \in E_M.
\end{equation*} 
The assignment to  $m \in M$ of
$(\alpha_m,e)\in \overline{\mathrm{Aut}(\mathrm{e}_M)}$
and to
$(\alpha,x)\in \overline{\mathrm{Aut}(\mathrm{e}_M)}$  of $x\in G$ yields a group extension
\begin{equation}
\begin{CD}
0
@>>>
M
@>{i_M}>>
\overline{\mathrm{Aut}(\mathrm{e}_M)}
@>{p_M}>>
G
@>>>
1.
\end{CD}
\label{extMM}
\end{equation}
The extension $\mathrm{e}_M$ being one of $N$-modules,
the $N$-action 
$n \mapsto \alpha_n\colon E_M\to E_M$ ($n \in N$)
on $E_M$ yields the injection 
\begin{equation}
s\colon N \longrightarrow \overline{\mathrm{Aut}(\mathrm{e}_M)},
\ 
n \longmapsto (\alpha_n,n),
\label{s}
\end{equation}
and
the composite of this injection with the projection to $G$
coincides with the injection $N \to G$ in \eqref{ge1}. 
By construction, the normalizer  $\mathrm{Aut}(\mathrm{e}_M)$ 
 of $s(N)$ in $\overline{\mathrm{Aut}(\mathrm{e}_M)}$
consists of the pairs
$(\alpha,x)\in \overline{\mathrm{Aut}(\mathrm{e}_M)}$ such that
\[
\alpha \alpha_n \alpha^{-1}=\alpha_{xnx^{-1}},\ n \in N.
\]
By \cite[Proposition 1.7]{MR597986}, $\mathrm{Aut}(\mathrm{e}_M)$
still maps onto $G$, and the surjection onto $G$ yields a group extension
\begin{equation}
\begin{CD}
0
@>>>
M^N
@>>>
\mathrm{Aut}(\mathrm{e}_M)
@>{\pi_{\mathrm{e}_M}}>>
G
@>>>
1.
\end{CD}
\label{extMMM}
\end{equation}
Passing to quotients, 
with the notation 
$\mathrm{Out}(\mathrm{e}_M)= \mathrm{Aut}(\mathrm{e}_M)/s(N)$,
we obtain the group extension
\begin{equation}
\begin{CD}
\widetilde {\mathrm{e}}_M
\colon
0
@>>>
M^N
@>>>
\mathrm{Out}(\mathrm{e}_M)
@>>>
Q
@>>>
1.
\end{CD}
\label{extMMMM}
\end{equation}
By \cite[Theorem 6]{MR597986}, the assignment 
of $\widetilde {\mathrm{e}}_M$
to $\mathrm{e}_M$
 yields a homomorphism
\begin{equation}
\Delta\colon 
\HH^1(N,M)^Q \longrightarrow 
\mathrm{Opext}(Q,M^N)
\label{Delta}
\end{equation}
and, via the canonical isomorphism $\mathrm{Opext}(Q,M^N)\cong \HH^2(Q,M^N)$,
the homomorphism $\Delta$ yields
a conceptual description of the differential
\eqref{d21}, viz.
$\mathrm{d_2}\colon 
\HH^1(N,M)^Q \longrightarrow 
\HH^2(Q,M^N)$.

\subsection{Relationship between the two previous descriptions}

Pulling back the group \eqref{sdfb} through the projection $G \to Q$
in \eqref{ge1}, we obtain the
group
\begin{equation}
\mathrm{Hom}_N(\mathbb Z G,M)\rtimes _TG
=\{(\varphi,x);\varphi|_{\mathrm I G}=
\alpha_x \mu \alpha_x^{-1}-\mu\}\subseteq 
\mathrm{Hom}_N(\mathbb Z G,M)\rtimes G . 
\label{pullb}
\end{equation}
Given $(\varphi,y) \in \mathrm{Hom}_N(\mathbb Z G,M)\rtimes G$,
define
\begin{align*}
\alpha_1&\colon \mathbb Z G \to E_M,\ \alpha_1(x)
=\varphi(yx)+\nu(yx),\ x \in \mathbb Z G,
\\
\alpha_2&\colon E_M \to E_M,\ \alpha_2(m)=ym .
\end{align*}
Since the left-hand square of 
\eqref{diag11} is a pushout diagram of abelian groups, $\alpha_1$ and 
$\alpha_2$ induce a unique automorphism 
\[
 \alpha_{\varphi,y}\colon E_M \longrightarrow E_M.
\]
By \cite[Lemma 5.1]{MR597986},
the assignment 
of  $\alpha_{\varphi,y}$
to 
$(\varphi,y)\in \mathrm{Hom}_N(\mathbb Z G,M)\rtimes _TG$
yields a morphism
\begin{equation}
\begin{CD}
0
@>>>
M^N
@>>>
 \mathrm{Hom}_N(\mathbb Z G,M)\rtimes _TG
@>>>
G
@>>>
1
\\
@.
@|
@V{\Phi}VV
@|
@.
\\
0
@>>>
M^N
@>>>
\mathrm{Aut}(\mathrm{e}_M)
@>>>
G
@>>>
1
\end{CD}
\label{diag9}
\end{equation}
of group extensions.
The resulting morphism
\begin{equation}
\begin{CD}
0
@>>>
M^N
@>>>
 \mathrm{Hom}_N(\mathbb Z G,M)\rtimes _TQ
@>>>
Q
@>>>
1
\\
@.
@|
@V{\widehat \Phi}VV
@|
@.
\\
0
@>>>
M^N
@>>>
\mathrm{Out}(\mathrm{e}_M)
@>>>
Q
@>>>
1
\end{CD}
\label{diag10}
\end{equation}
of group extensions identifies
$\mathrm{Out}(\mathrm{e}_M)$ with  $\mathrm{Hom}_N(\mathbb Z G,M)\rtimes _TQ$.
In view of 
\cite[Proposition 3.1]{MR641328},
diagram \eqref{diag10}  entails a proof of Theorem 6  
in \cite{MR597986}.

\subsection{Construction in \cite{MR2959787}}
\label{section5}

The construction of 
\eqref{tr}
in \cite[Section 5]{MR2959787}
comes down to this:
Consider the split extension
\begin{equation}
\begin{CD}
0
@>>>
M
@>>>
M\rtimes G
@>>>
G
@>>>
1.
\end{CD}
\label{extsplit}
\end{equation}
A section $s \colon N \to M\rtimes G$
whose composite with the projection to $G$ coincides with
the injection $N \to G$ in 
$\mathrm{e}_G\colon
1
\to
N
\stackrel{i^N}
\to
G
\stackrel{\pi}
\to
Q
\to
1
$
is determined by its component
$\varphi \colon N \to M$, and $\varphi$ is a
1-cocycle or, equivalently, derivation, 
on $N$ with values in $M$
that represents a $Q$-invariant cohomology class
$[\varphi]\in \HH^1(N,M)^Q$.
Given such a section $s$ or, equivalently, 1-cocycle $\varphi$,
the normalizer $N_{M\rtimes G}(sN)$
of $sN$ in $M\rtimes G$ maps onto $G$,  the sequence
\begin{equation}
\begin{CD}
0
@>>>
M^N
@>>>
N_{M\rtimes G}(sN)
@>>>
G
@>>>
1
\end{CD}
\label{exact4}
\end{equation}
is exact, and
so is
\begin{equation}
\begin{CD}
\widetilde{\mathrm e}\colon
0
@>>>
M^N
@>>>
N_{M\rtimes G}(sN)/sN
@>>>
Q
@>>>
1.
\end{CD}
\label{etilde}
\end{equation}
In
 \cite{MR2959787}
the map \eqref{tr}
is essentially 
defined by the assignment
to an $M$-valued 1-cocycle $\varphi$ on $N$
representing a $Q$-invariant cohomology class in $\HH^1(N,M)$
of the class of the extension
$\widetilde{\mathrm e}$. 
The construction in  \cite{MR2959787}
is actually couched in what is referred to as \lq\lq partial semi-direct
complements\rq\rq.
The 1-cocycle $\varphi$ is lurking behind the choice of 
the subgroup of $M\rtimes G$  denoted  in  \cite{MR2959787} by $H$,
indeed, it is the component into $M$ of the resulting
section $s\colon N \to M\rtimes G$ such that $H=sN$.

\subsection{Comparison of the constructions in \cite{MR597986} and
\cite{MR2959787} and proof of Theorem \ref{one}}
\label{compar}
We maintain the notation established before.
Thus $\mathrm{e}_M 
\colon
0
\to
M
\stackrel{i_M}
\to
E_M
\stackrel{p_M}
\to
\mathbb Z
\to
0
$
is an extension of $N$-modules of the kind
\eqref{extM}
that represents a class
\[
[\mathrm{e}_M]\in 
\mathrm{Ext}_N^1(\mathbb Z,M)^Q\cong \HH^1(N,M)^Q,
\]
and $s\colon N \to M\rtimes G$ is a section as in Subsection \ref{section5}
above.

As an abelian group, $E_M$ decomposes as a direct sum $M \oplus \mathbb Z$.
Such a decomposition induces an obvious $G$-module structure on  $E_M$
that restricts to the $G$-module structure on $M$.
This $G$-module structure induces a section
$\sigma\colon G \to \overline{\mathrm{Aut}(\mathrm{e}_M)}$
for the extension \eqref{extMM} and hence an isomorphism
$\overline{\mathrm{Aut}(\mathrm{e}_M)} \to M \rtimes G$.
This isomorphism induces a morphism
\begin{equation}
\begin{CD}
\widetilde{\mathrm e}\colon
0
@>>>
M^N
@>>>
N_{M\rtimes G}(sN)/sN
@>>>
Q
@>>>
1
\\
@.
@|
@VVV
@|
@.
\\
\widetilde {\mathrm{e}}_M\colon
0
@>>>
M^N
@>>>
\mathrm{Out}(\mathrm{e}_M)
@>>>
Q
@>>>
1
\end{CD}
\label{diag20}
\end{equation}
of group extensions.
Hence, cf. the introductory phrase
of Subsection \ref{hueb} above,
the construction in  
\cite[Section 5]{MR2959787} that underlies the map
\eqref{tr} is exactly the same as that in
\cite{MR597986} which underlies the map \eqref{Delta}.
The only difference is that the interpretation of
the naive semi-direct product $M \rtimes G$
as the group $\overline{\mathrm{Aut}(\mathrm{e}_M)}$,
not noticed in \cite{MR2959787},
enabled us to show in \cite{MR597986} that \eqref{Delta}
actually yields the differential
\eqref{d21}.
In particular, the map \eqref{tr}
coincides with the differential \eqref{d21}.
This establishes Theorem \ref{one}.
Furthermore, the claim
in \cite{MR2995036}
to the effect that the construction of the map \eqref{tr}
in \cite{MR2959787} was much more elementary than the construction
of \eqref{Delta} in \cite{MR597986} 
is unfounded.

\section{$\HH^2(G,M)_1 \longrightarrow \HH^1(Q,\HH^1(N,M))$} 
\label{d23}

We note first that, the actions of the various groups on
a module being written as juxtaposition, the standard
$G$-module structure on $\mathrm{Der}(N,M)$
is given by the familiar formula
\[
(yd)(n) = y(d(y^{-1}n y)),\ d \in \mathrm{Der}(N,M),\ y \in G,\ n \in N;
\]
this $G$-module structure 
induces,
through the projection from $E$ to $G$,
 an  $E$-module structure
on $\mathrm{Der}(N,M)$, and this structure, in turn, 
induces the standard
 $Q$-module structure
on $\HH^1(N,M)$.

\subsection{Construction of the map
$\rho\colon 
\HH^2(G,M)_1 \to
\HH^1(Q,\HH^1(N,M))
$
in  \cite{MR2959787}}
Consider a group extension
\begin{equation}
\begin{CD}
\eG\colon 
0
@>>>
M
@>i>>
E
@>p>>
G
@>>>
1.
\end{CD}
\label{ext5}
\end{equation}
Suppose that the restriction 
of $\eG$
to $N$ splits, and let
$s\colon N \to E$ be a homomorphism such that
$p\circ s$ coincides with the injection of $N$ into $G$.
Given $x \in E$, we denote its image $p(x)\in G$ by $[x]$. 
The assignment to $x \in E$ of
\begin{equation}
\widehat \delta_x\colon N \longrightarrow M,\ 
i (\widehat \delta_x(n))=xs([x^{-1}]n[x])x^{-1}s(n^{-1})\in i(M) \subseteq E,\ n \in N,
\label{der3}
\end{equation}
yields a derivation 
$\widehat \delta_{\eG}\colon E \to \mathrm{Der}(N,M)$.
Given $x \in E_N=p^{-1}(N)$, the member
$y=(s[x])x^{-1}$ of $E_N$ lies in the image $i(M)$ of $M$ in $E_N$, and
\[
i (\widehat \delta_x(n)) =y^{-1} s(n) y s(n^{-1}) \in i(M) \subseteq E_N,\  n \in N,
\]
that is, $\widehat \delta_x$ is an inner derivation.
Consequently
the derivation $\widehat \delta_{\eG}$, in turn, passes to a derivation
\begin{equation}
\delta_{\eG}\colon Q \longrightarrow \HH^1(N,M) 
\label{der1}
\end{equation}
and, in  
\cite[Section 6]{MR2959787}, the map
$\rho\colon 
\HH^2(G,M)_1 \to
\HH^1(Q,\HH^1(N,M))$ is defined by the
assignment of $\delta_{\eG}$ to $\eG$.
While the derivations
$\widehat \delta_{\eG}$
and
$\delta_{\eG}$ depend on the choice of the homomorphism ${s\colon N \to E}$,
the class $[\delta_{\eG}]\in \HH^1(Q,\HH^1(N,M))$ is independent of that 
choice.

\subsection{Crossed pairs}
At the referee's request,
for intelligibility,
we will now briefly reproduce 
some material from \cite{MR597986}.
To this end, we write the $G$-action on $M$ as
\begin{equation*}
G \times M \longrightarrow M,\ (x,y) \longmapsto  {}^xy,\ x \in G,\ y \in M.
\end{equation*}
Let
\begin{equation*}
\mathrm e \colon 0 \longrightarrow M \longrightarrow \EN \stackrel{\piN} \longrightarrow N \longrightarrow 1
\end{equation*}
be a group extension
whose class $[\mathrm e ] \in \mathrm H^2 (N,M)$ is fixed under the standard $Q$-action
on $\mathrm H^2(N,M)$. 
Given $x\in G$, we write
\begin{equation*}
\ell_x(y) = {}^xy, \ y \in M,\quad 
i_x(n) = xnx^{-1},\ n \in N.
\end{equation*}
Write $\Aut_G(\mathrm e) $ for the subgroup
of $\Aut (\EN) \times G$ that consists of pairs $(\kappa , x)$ which
make the diagram
\begin{equation*}
\begin{CD}
0
@>>> 
M 
@>>> 
\EN 
@>>> 
N 
@>>> 
1
\\
@.
@V{\ell_x}VV
@V{\kappa}VV
@V{i_x}VV
@.
\\
0
@>>> 
M 
@>>> 
\EN 
@>>> 
N 
@>>> 
1
\end{CD}
\end{equation*}
commutative.

The homomorphism 
\begin{equation*}
\beta\colon \EN \longrightarrow \Aut_G(\mathrm e ),\ 
\beta(y)=(i_y,i^N(\piN(y))), \ y \in \EN,
\end{equation*}
together with the obvious action
of $\Aut_G(\mathrm e )$ on $\EN$, yields a crossed module
$(\EN,\Aut_G(\mathrm e ),\beta)$ whence, in particular,
$\beta(\EN)$ is a normal subgroup of $\Aut_G(\mathrm e )$;
with the notation $\Out_G(\mathrm e)$
for the cokernel of $\beta$, 
the resulting crossed 2-fold extension has the form
\begin{equation}
\begin{CD}
0 @>>>  
M^N
 @>>> \EN @>{\beta}>> \Aut_G(\mathrm e ) @>>> 
\Out_G(\mathrm e) @>>> 1.
\end{CD}
\label{res1}
\end{equation}
The association
\begin{equation}
\mathrm{Der}(N,M) \ni d\longmapsto (\kappa_d,1),\ \kappa_d(y)= (d\piN(y))y,\ 
y \in \EN,
\label{vartheta}
\end{equation} 
yields an injective homomorphism
$\vartheta\colon \mathrm{Der}(N,M) \longrightarrow \Aut_G(\mathrm e)$;
this homomorphism and the obvious map $\Aut_G(\mathrm e) \to G$
yield a group extension
\begin{equation}
0 \longrightarrow \Der (N,M) 
\stackrel{\vartheta}
\longrightarrow \Aut_G(\mathrm e ) \longrightarrow G \longrightarrow 1,
\label{ext9}
\end{equation}
the map $\Aut_G(\mathrm e) \to G$ being surjective,
since the  class $[\mathrm e ] \in \mathrm H^2 (N,M)$ is supposed to be
fixed under $Q$.
Further,
the formula
$(\zeta(m))(n)=m({}^nm)^{-1}$, as $m$ ranges over $M$ and $n$ over $N$,
defines a homomorphism 
$\zeta\colon M \to \mathrm{Der}(N,M)$. 
With these preparations out of the way, the various groups under discussion 
fit into the commutative
diagram
\begin{equation}
\begin{CD}
@. 0 @. 0 @.@.
\\
@.
@VVV
@VVV
@.
@.
\\
@. 
M^N 
@= 
M^N
@. 
1
@. 
\\
@.
@VVV
@VVV
@VVV
@.
\\
0
@>>> 
M 
@>>> 
\EN 
@>{\piN}>> 
N 
@>>> 
1
\\
@.
@V{\zeta}VV
@V{\beta}VV
@V{i^N}VV
@.
\\
0 @>>> \Der (N,M) @>{\vartheta}>>       \Aut_G(\mathrm e ) @>>> G @>>> 1
\\
@.
@VVV
@VVV
@VVV
@.
\\
0 @>>> \mathrm H^1(N,M) @>>> \Out_G(\mathrm e ) @>>> Q @>>> 1
\\
@.
@VVV
@VVV
@VVV
@.
\\
@. 1 @. 1 @.1 @.
\end{CD}
\label{diag1}
\end{equation}
with exact rows and columns.
We shall use the notation
\begin{equation}
\begin{CD}
\overline {\mathrm e} \colon 
0 @>>> \mathrm H^1(N,M) @>>> \Out_G(\mathrm e ) @>>> Q @>>> 1
\end{CD}
\label{ext8}
\end{equation}
for the bottom row extension of \eqref{diag1}.
This extension is the cokernel, in the category
of group extensions with abelian kernel,
of the morphism
$(\zeta,\beta,i)$ of group extensions.

Suppose now that the
group extension $\overline{\mathrm e}$
splits; we then say that the group extension
$\mathrm e\colon 0 \to M \to \EN \stackrel{\piN} \to N \to 1$ 
{\em admits a crossed pair structure with respect to the group extension\/}
\eqref{ge1},
and we refer to 
a section $\psi\colon Q \to  \Out_G(\mathrm e )$
of $\overline{\mathrm e}$ as
a {\em crossed pair structure
on the  group extension $\mathrm e$ 
with respect to the group extension\/}
\eqref{ge1}.
By definition, a {\em crossed pair\/}
$(\mathrm e,\psi)$ {\em with respect to the group extension\/}
\eqref{ge1} {\em and the\/} $G$-{\em module\/} $M$
consists of a group extension
$\mathrm e\colon 0 \to M \to \EN \to N \to 1$
whose class $[\mathrm e]\in \mathrm H^2(N,M)$ is fixed under $Q$
such that  the associated extension $\overline{\mathrm e}$ splits,
together with a section 
$\psi\colon Q \to  \Out_G(\mathrm e )$
of $\overline{\mathrm e}$ \cite[p.~152]{MR597986}.

\begin{Remark} {\rm By 
\cite[Theorem~1]{MR641328},
the association $\mathrm e \mapsto \overline {\mathrm e}$
yields a conceptual 
description of the differential 
$d_2 \colon \mathrm E^{0,2}_2(M) \to \mathrm E^{2,1}_2(M)$
of 
the 
spectral sequence 
$(\mathrm E^{p,q}_r(M),d_r)$
associated to the group extension \eqref{ge1} and the $G$-module $M$.}
\end{Remark}

Under a kind of generalized Baer sum, 
suitable
classes of crossed pairs with respect to 
\eqref{ge1} and the $G$-module $M$
 constitute an abelian group 
$\mathrm{Xpext}(G,N; M)$ \cite[Theorem~1]{MR597986}. 
Moreover, 
cf. 
 \cite[Theorem~2]{MR597986},
suitably defined homomorphisms
\begin{equation*}
j\colon \mathrm H^2(G,M)\longrightarrow
\mathrm{Xpext} (G,N;M),\ 
\Delta\colon \mathrm{Xpext} (G,N;M)
\longrightarrow \mathrm H^3(Q,M^N)
\end{equation*}
yield an extension of the classical five term exact sequence to an
eight term exact sequence of the kind
\eqref{exact3}.

To recall the construction of $j$,
consider a group extension 
\begin{equation*}
\eG \colon 0 \longrightarrow M 
\stackrel{i}
\longrightarrow E 
\stackrel{p}
\longrightarrow 
G \longrightarrow 1;
\end{equation*}
the restriction 
of $\eG$
to $N$ has the form 
$\eMN \colon 0 \to M \to E_N \to N \to 1$,
and conjugation in $E$ induces a crossed pair structure
$\psi_E\colon Q \to  \Out_G(\eMN )$ 
on $\eMN$ 
(section for \eqref{ext8})
with respect to \eqref{ge1}.
The map $j$ is given by the assignment to a group extension of the kind
$\eG$ of the crossed pair
\begin{equation*}
\left(\eMN\colon 0 \to M \to E_N  \to N \to 1, 
\ \psi_E\colon Q \to  \Out_G(\eMN )\right) 
\end{equation*}
with respect to \eqref{ge1}.

To recall the construction of $\Delta$,
given a crossed pair
\begin{equation*}
\left(\mathrm e\colon 0 \to M \to \EN  \to N \to 1, 
\ \psi\colon Q \to  \Out_G(\mathrm e )\right)
\end{equation*}
with respect to the group extension \eqref{ge1} and the $G$-module $M$,
let $B^\psi $ denote the fiber product group
$\Aut_G (\mathrm e ) \times_{\Out_G (\mathrm e)} Q$
with respect to the crossed pair structure map
$\psi\colon Q \to \Out_G (\mathrm e)$
and, furthermore, let
$\partial^\psi\colon
\EN \to B^\psi$ denote
the obvious homomorphism; together with the obvious
action of $B^\psi$ on $\EN$ induced by the canonical homomorphism
 $B^\psi \to \Aut_G (\mathrm e )$,
the exact sequence
\begin{equation}
\mathrm e_\psi \colon 0 \longrightarrow M^N \longrightarrow \EN \stackrel{\partial^\psi}{\longrightarrow} B^\psi \longrightarrow
Q \longrightarrow 1
\label{epsi}
\end{equation}
is a crossed 2-fold extension and hence represents a class
in $\mathrm H^3(Q,M^N)$.
We shall refer to
$\mathrm e_\psi$ {\em as the crossed\/} 2-{\em fold extension associated to
the crossed pair\/} $(\mathrm e,\psi)$.
The homomorphism $\Delta\colon 
\mathrm{Xpext} (G,N;M) \to\mathrm H^3(Q,M^N)$
is given by the assignment to
a crossed pair
$(\mathrm e,\psi)$ 
with respect to the group extension \eqref{ge1} and the $G$-module $M$
of its associated crossed 2-fold extension
$\mathrm e_\psi$.

Next we recall the injection $\alpha\colon
\HH^1(Q,\HH^1(N,M)) \longrightarrow \mathrm{Xpext} (G,N;M)$
constructed in 
\cite[Subsection 1.3]{MR597986}. To this end, let $\mathrm e_s$ denote the 
split extension
\begin{equation*}
\mathrm e_s \colon 0 \longrightarrow M \longrightarrow M\rtimes N 
\longrightarrow 
N \longrightarrow 1.
\end{equation*}
The obvious action of $G$ on
$M\rtimes N$ induces, on
$\mathrm e_s$,
 a canonical crossed pair structure
${\psi_s\colon Q \to \Out_G(\mathrm e_s)}$
with respect to \eqref{ge1}.
Given a derivation $\delta\colon Q \to \HH^1(N,M)$, define the
crossed pair structure 
\begin{equation*}
\psi_{\delta} \colon Q \longrightarrow \Out_G(\mathrm e_s)
\end{equation*}
on $\mathrm e_s$ with respect to \eqref{ge1}
by 
\begin{equation}
\psi_{\delta}(x) =\delta(x) \psi_s(x), \ x \in Q;
\label{cpd}
\end{equation}
here we identify in notation the member
$\delta(x)$ of $\HH^1(N,M)$ with its image in
$\Out_G(\mathrm e)$ under the injection
$\HH^1(N,M) \to \Out_G(\mathrm e)$ in \eqref{ext8}.
The assignment to a derivation ${\delta\colon Q \to \HH^1(N,M)}$
of the crossed pair $(\mathrm e_s, \psi_{\delta})$
with respect to \eqref{ge1} and the $G$-module $M$
induces a homomorphism 
\begin{equation*}
\alpha\colon
\HH^1(Q,\HH^1(N,M)) \longrightarrow \mathrm{Xpext} (G,N;M).
\end{equation*}
Let 
\[{\xi\colon \mathrm{Xpext} (G,N;M) \longrightarrow \HH^0(Q,\HH^2(N,M))}
\]
denote the canonical forgetful map.
By
\cite[Theorem 3]{MR597986}, the sequence
\begin{equation}
0
\to
\HH^1(Q,\HH^1(N,M))
\stackrel{\alpha}
\to
\mathrm{Xpext} (G,N;M) 
\stackrel{\xi}
\to 
\HH^0(Q,\HH^2(N,M))
\stackrel{d_2}
\to
\HH^2(Q,\HH^1(N,M))
\label{exact2}
\end{equation}
is exact and, furthermore, natural in the data.
The exactness of \eqref{exact2} at $\mathrm{Xpext} (G,N;M)$
entails that, via $\alpha$, the group 
$\HH^1(Q,\HH^1(N,M))$ acts faithfully and transitively on the
classes of crossed pair structures 
$\psi\colon Q \to \Out_G(\mathrm e_s)$
on the split extension 
\[
\mathrm e_s \colon 0 \longrightarrow M \longrightarrow M \rtimes N \longrightarrow N \longrightarrow 1.
\]
Somewhat more explicitly, this action is characterized by the bijective 
correspondence between classes of crossed pair structures
on $\mathrm e_s$ and cohomology classes of
derivations $Q \to \HH^1(N,M)$
induced by  \eqref{cpd} above.

\subsection{Proof of Theorem \ref{two}}
Below we will use the notation 
{\rm \lq\lq can\rq\rq}\ 
for 
any map  that is canonically defined.

\begin{prop}
\label{prop2}
The diagram
\begin{equation}
\begin{CD}
\HH^2(G,M)_1 @>{\mathrm{can}}>> \HH^1(Q,\HH^1(N,M))
@>{d_2}>>
\HH^3(Q,M^N)
\\
@V{\mathrm{can}}VV
@V{\alpha}VV
@|
\\
\HH^2(G,M) @>{j}>>  \mathrm{Xpext}(G,N;M)
@>{\Delta}>>
\HH^3(Q,M^N)
\\
@V{\mathrm{can}}VV
@V{\xi}VV
@V{\mathrm{can}}VV
\\
\mathrm E_{\infty}^{0,2}(M)
@>{\mathrm{can}}>>
\mathrm E_{3}^{0,2}(M)
@>{d_3}>>
\mathrm E_{3}^{3,0}(M)
\end{CD}
\label{diag22}
\end{equation}
is commutative with exact rows.
\end{prop}
\begin{proof} 
This is a consequence of the commutativity of
diagram \cite[(1.10)]{MR597986}.
\end{proof}

\begin{prop}
\label{prop3}
The diagram
\begin{equation}
\begin{CD}
\HH^2(G,M)_1 @>{\rho}>> \HH^1(Q,\HH^1(N,M))
@>{d_2}>>
\HH^3(Q,M^N)
\\
@V{\mathrm{can}}VV
@V{\alpha}VV
@|
\\
\HH^2(G,M) @>{j}>>  \mathrm{Xpext}(G,N;M)
@>{\Delta}>>
\HH^3(Q,M^N)
\end{CD}
\label{diag21}
\end{equation}
is commutative
with exact rows.
\end{prop}

\begin{proof}
The commutativity of the right-hand square
is implied by the commutativity of
diagram \eqref{diag22}. 

We will now show that
the left-hand square in \eqref{diag21}
is commutative.
To this end,
consider a group extension 
\begin{equation*}
\eG \colon 0 \longrightarrow M \stackrel{i}\longrightarrow E 
\stackrel{p}\longrightarrow 
G \longrightarrow 1
\end{equation*}
that represents a member of $\HH^2(G,M)_1$. 
Then the restriction 
$\eMN \colon 0 \to M \to E_N \to N \to 1$
of $\eG$
to $N$ splits; choose a section $s \colon N \to E_N$,
and write  
$E_N$ as $M \rtimes (sN)$.
By the construction of $j$,
the group $M \rtimes (sN)$ being normal in $E$,
conjugation in $E$ induces a crossed pair structure
$\psi_E\colon Q \to  \Out_G(\eMN )$ 
on $\eMN$ (section for \eqref{ext8}) with respect to \eqref{ge1},
and the crossed pair $(\eMN,\psi_E)$
with respect to \eqref{ge1}
represents 
the
image $j[\eG] \in  \mathrm{Xpext}(G,N;M)$
of the class  $[\eG] \in \HH^2(G,M)_1$.
By construction, $\xi(j[\eG]) =0$ whence 
 \[
j[\eG] \in  \mathrm{ker}(\xi)
=\alpha(\HH^1(Q,\HH^1(N,M)))
\subseteq \mathrm{Xpext}(G,N;M).
\]

The canonical action of $G$ on $M \rtimes (sN)$ 
induces a section 
$\Psi_s\colon G \to \Aut_G(\eMN)$
for \eqref{ext9}. 
On the other hand, conjugation in $E$ induces likewise a section
$\Psi_{\eG}\colon E \to \Aut_G(\eMN)$
for \eqref{ext9}, and this section, in turn,
induces the crossed pair structure
$\psi_E\colon Q \to  \Out_G(\eMN )$ 
on $\eMN$ with respect to \eqref{ge1}.
Recall  the injection $\vartheta\colon \mathrm{Der}(N,M) \to \Aut_G(\eMN)$
characterized by \eqref{vartheta}.

\begin{prop}
\label{prop}
The derivation $\widehat \delta_{\eG}\colon E \to \mathrm{Der}(N,M)$
satisfies the identity
\begin{equation}
\Psi_{\eG} (x)=
\vartheta(\widehat \delta_{\eG}(x))\Psi_s([x]),\ x \in E.
\label{ident3}
\end{equation}
\end{prop}

\begin{proof}
Define
the derivation $\widetilde \delta_{\eG}\colon E \to \mathrm{Der}(N,M)$
by the identity
\begin{equation}
\Psi_{\eG} (x)=
\vartheta(\widetilde \delta_{\eG}(x))\Psi_s([x]),\ x \in E.
\label{ident33}
\end{equation}
We shall show that
$\widehat \delta_{\eG}= \widetilde \delta_{\eG}$.

By construction, an explicit expression for 
$\Psi_s\colon G \to \Aut_G(\eMN)$ is given by
\begin{equation*}
(\Psi_s(x))(m,s(n))=({}^xm,s(xnx^{-1})),\ x\in G,\ m \in M,\ n \in N.
\end{equation*}
Let $x \in E$, $m \in M$, and $n \in N$.
Since $i({}^{[x]}m)=xi(m)x^{-1}$, 
in view of the identity \eqref{vartheta}
characterizing 
the injection $\vartheta\colon \mathrm{Der}(N,M) \to \Aut_G(\eMN)$,
\begin{align*}
\vartheta(\widetilde \delta_{\eG}(x))\Psi_s([x])
(i(m)s(n))&=\vartheta(\widetilde \delta_{\eG}(x))
(xi(m)x^{-1}
s([x]n[x]^{-1}))
\\&=
i(\widetilde \delta_{\eG}(x)([x]n[x]^{-1}))
xi(m)x^{-1}
s([x]n[x]^{-1}).
\end{align*}
Since, by construction,
\begin{equation}
(\Psi_{\eG} (x))(i(m)s(n))=i({}^{[x]}m)xs(n)x^{-1},
\label{ident6}
\end{equation}
from \eqref{ident33}, we deduce the identity
\begin{equation}
i({}^{[x]}m)xs(n)x^{-1}
=
i(\widetilde \delta_{\eG}(x)([x]n[x]^{-1}))
xi(m)x^{-1}
s([x]n[x]^{-1}). 
\label{ident5}
\end{equation}
To simplify the notation, we will now write
$\widetilde \delta =\widetilde \delta_{\eG}$
and $\widetilde \delta_x=\widetilde \delta(x)$.
Then the identity \eqref{ident5}
takes the form
\begin{equation}
i({}^{[x]}m)xs(n)x^{-1}
=
i(\widetilde \delta_x([x]n[x]^{-1}))
xi(m)x^{-1}
s([x]n[x]^{-1}) 
\label{ident7}
\end{equation}
or, since $M$ is abelian,
\begin{equation}
i({}^{[x]}m)xs(n)x^{-1}
=
i({}^{[x]}m)
i(\widetilde \delta_x([x]n[x]^{-1}))
s([x]n[x]^{-1}), 
\label{ident8}
\end{equation}
and thence
\begin{equation}
xs(n)x^{-1}
=
i(\widetilde \delta_x([x]n[x]^{-1}))
s([x]n[x]^{-1}), 
\label{ident9}
\end{equation}
whence
\begin{equation}
i(\widetilde \delta_x([x]n[x]^{-1}))
=
xs(n)x^{-1}
s([x]n^{-1}[x]^{-1}), 
\label{ident10}
\end{equation}
or, equivalently,
\begin{equation}
i(\widetilde \delta_x(n))
=
x
(s([x]^{-1}n[x]))
x^{-1}
s(n^{-1}). 
\label{ident11}
\end{equation}
In view of the identity \eqref{der3}
characterizing
$\widehat \delta_{\eG}$, we conclude that
$\widehat \delta_{\eG}= \widetilde \delta_{\eG}$.
\end{proof}

We now complete the proof of Proposition \ref{prop3}.
The commutativity
of the left-hand square in \eqref{diag21}
is a consequence of Proposition \ref{prop}.
Indeed,
once a choice of section $s$ for 
the restriction
$\eMN \colon 0 \to M \to E_N \to N \to 1$
of $\eG\colon 
 0 \to M \stackrel{i}\to E 
\stackrel{p}\to 
G \to 1
$
to $N$
has been made,
the derivation 
$\delta_{\eG}\colon Q \to \HH^1(N,M)$
represents the image $\rho([\eG])\in \HH^1(Q,\HH^1(N,M))$
of $[\eG] \in \HH^2(G,M)_1$. Now the section 
$\Psi_{\eG}\colon E \to \Aut_G(\eMN)$ 
for \eqref{ext9}
induces the
crossed pair structure $\psi_E\colon Q \to  \Out_G(\eMN)$ 
on $\eMN$
that yields the image
$j([\eG]) \in \mathrm{Xpext}(G,N;M)$.
On the other hand,
the derivation $\widehat \delta_{\eG}\colon E \to \mathrm{Der}(N,M)$
induces the derivation ${\delta_{\eG}\colon Q \to \HH^1(N,M)}$.
Hence the identity \eqref{ident3}
entails that, by its very construction, the 
injective homomorphism
${\alpha\colon
\HH^1(Q,\HH^1(N,M)) \to \mathrm{Xpext} (G,N;M)}$
recovers 
the class of the crossed pair structure 
$\psi_E\colon Q \to \Out_G(\eMN)$
on $\eMN$
from 
$\rho([\eG])=[\delta_{\eG}]\in \HH^1(Q,\HH^1(N,M))$.
\end{proof}

The homomorphism $\alpha\colon
\HH^1(Q,\HH^1(N,M)) \longrightarrow \mathrm{Xpext} (G,N;M)$ 
is injective whence,
by the commutativity of the left-hand square
in \eqref{diag21},
$\rho\colon \HH^2(G,M)_1\to \HH^1(Q,\HH^1(N,M))$
arises as the restriction of $j\colon \HH^2(G,M)\to \mathrm{Xpext} (G,N;M)$
to $\HH^2(G,M)_1$,
and the exact sequence in \cite{MR2959787}
can be deduced from the exact sequence
\eqref{exact3} constructed in  \cite{MR597986}.
The commutativity of \eqref{diag21} 
and the exactness of the bottom row
of that diagram
entail that the homomorphism
$\rho\colon \HH^2(G,M)_1\to \HH^1(Q,\HH^1(N,M))$ 
coincides with the canonical homomorphism
\begin{equation}
\mathrm{can}\colon
\HH^2(G,M)_1
\longrightarrow 
\mathrm{E}_{\infty}^{1,1}(M)
\longrightarrow
\mathrm{E}_{2}^{1,1}(M)
\cong
\HH^1(Q,\HH^1(N,M))
\label{can}
\end{equation}
that comes
from the spectral sequence $(\mathrm E^{p,q}_r(M),d_r)$.
This establishes Theorem \ref{two}.
In particular,
$\rho\colon \HH^2(G,M)_1\to \HH^1(Q,\HH^1(N,M))$ 
maps $\HH^2(G,M)_1$ onto 
\[
\mathrm E^{1,1}_{\infty}=\ker (d_2\colon
 \HH^1(Q,\HH^1(N,M))
\to
\HH^3(Q,M^N)).
\] 

\begin{Remark}{\rm
A derivation or 1-cocycle of the kind \eqref{der1}
or of the kind $\delta$ in \eqref{cpd}
can be seen as a kind of {\em Eilenberg difference cocycle\/}
\cite[p.~236]{MR0001351}.
The identity \eqref{der3} 
spells out $\widehat \delta_x$ 
as the \lq\lq difference\rq\rq\ of the two sections
${}^x s$ and $s$ from $N$ to $E$ for the restriction of the
group extension $\eG\colon 0 \to M \to E \to G \to 1$ to $N$.
The key step in the proof of Theorem \ref{two}
consists in making the requisite difference cocycle explicit via the identity
\eqref{ident3}; this identity displays 
the 1-cocycle $\widehat \delta_{\eG}\colon E \to \Der(N,M)$ as the \lq\lq difference\rq\rq\ 
of the two sections
$\Psi_{\eG}$
and $\Psi_s$ for the corresponding extension \eqref{ext9}.}
\end{Remark}

\addcontentsline{toc}{section}{References}
\bibliographystyle{alpha}

\newcommand{\etalchar}[1]{$^{#1}$}
\def\cprime{$'$} \def\cprime{$'$} \def\dbar{\leavevmode\hbox to 0pt{\hskip.2ex
  \accent"16\hss}d} \def\cprime{$'$} \def\cprime{$'$} \def\cprime{$'$}
  \def\cprime{$'$} \def\cprime{$'$} \def\Dbar{\leavevmode\lower.6ex\hbox to
  0pt{\hskip-.23ex \accent"16\hss}D}
  \def\cftil#1{\ifmmode\setbox7\hbox{$\accent"5E#1$}\else
  \setbox7\hbox{\accent"5E#1}\penalty 10000\relax\fi\raise 1\ht7
  \hbox{\lower1.15ex\hbox to 1\wd7{\hss\accent"7E\hss}}\penalty 10000
  \hskip-1\wd7\penalty 10000\box7}
  \def\cfudot#1{\ifmmode\setbox7\hbox{$\accent"5E#1$}\else
  \setbox7\hbox{\accent"5E#1}\penalty 10000\relax\fi\raise 1\ht7
  \hbox{\raise.1ex\hbox to 1\wd7{\hss.\hss}}\penalty 10000 \hskip-1\wd7\penalty
  10000\box7} \def\polhk#1{\setbox0=\hbox{#1}{\ooalign{\hidewidth
  \lower1.5ex\hbox{`}\hidewidth\crcr\unhbox0}}}
  \def\polhk#1{\setbox0=\hbox{#1}{\ooalign{\hidewidth
  \lower1.5ex\hbox{`}\hidewidth\crcr\unhbox0}}}
  \def\polhk#1{\setbox0=\hbox{#1}{\ooalign{\hidewidth
  \lower1.5ex\hbox{`}\hidewidth\crcr\unhbox0}}}

\end{document}